\documentclass[reqno]{amsart}
\usepackage {amsmath,amssymb}
\usepackage {color,graphicx}
\usepackage [latin1]{inputenc}
\usepackage {times}
\usepackage {float}
\usepackage[dvips]{hyperref}
\usepackage{amscd}
\usepackage[all]{xy}
\usepackage{color}

\DeclareMathOperator {\id}{id}

\DeclareMathOperator {\Div}{Div}

\newcommand {\RR}{{\mathbb R}}

\newcommand {\arxiv}[3]{#1, \emph{#2}, preprint \href{http://arxiv.org/abs/#3}{arxiv:#3}.}

\newcommand {\arxivjournal}[4]{#1, \emph{#2}, #3; also at \href{http://arxiv.org/abs/#4}{arxiv:#4}.}

\newtheorem {theorem}{Theorem}[section]
\newtheorem {proposition}[theorem]{Proposition}
\newtheorem {lemma}[theorem]{Lemma}

\newtheorem {corollary}[theorem]{Corollary}
\theoremstyle {definition}
\newtheorem {definition}[theorem]{Definition}
\newtheorem {example}[theorem]{Example}

\newtheorem {notation}[theorem]{Notation}
\theoremstyle {remark}
\newtheorem {remark}[theorem]{Remark}

\begin {document}

\title {Tropical intersection products on smooth varieties}
\author {Lars Allermann}
\address {Lars Allermann, Fachbereich Mathematik, TU Kaiserslautern, Postfach 3049, 67653 Kaiserslautern, Germany}
\email {allerman@mathematik.uni-kl.de}

\begin{abstract}
In analogy to \cite[chapter 9]{AR07} we define an intersection
product of tropical cycles on tropical linear spaces $L^n_k$, i.e.
on tropical fans of the type $\max \{ 0, x_1,\ldots,x_n \}^{n-k}
\cdot \RR^n$. Afterwards we use this result to obtain an
intersection product of cycles on every smooth tropical variety,
i.e. on every tropical variety that arises from gluing such
tropical linear spaces. In contrast to classical algebraic
geometry these products always yield well-defined cycles, not
cycle classes only. Using these intersection products we are able
to define the pull-back of a tropical cycle along a morphism
between smooth tropical varieties. In the present article we stick
to the definitions, notions and concepts introduced in
\cite{AR07}.
\end{abstract}

\maketitle

\section{Intersection products on tropical linear spaces} \label{sec-tropicallinearspaces}

In this section we will give a proof that tropical linear spaces
$L^n_k$ admit an intersection product. Therefore we show at first
that the diagonal in the Cartesian product $L^n_k \times L^n_k$ of
such a linear space with itself is a sum of products of Cartier
divisors. Given two cycles $C$ and $D$ we can then intersect the
diagonal with $C \times D$ and define the product $C \cdot D$ to
be the projection thereof.

Throughout the section $e_1,\dots,e_n$ will always be the standard
basis vectors in $\RR^n$ and $e_0 := -e_1- \ldots -e_n$.

We begin the section with our basic definitions:

\begin{definition}[Tropical linear spaces] \label{def-Lnk}
  For $I \subsetneq
  \{0,1,\dots,n\}$ let $\sigma_I$ be the cone generated by the
  vectors $e_i$, $i \in I$. We denote by $L^n_k$ the tropical
  fan consisting of all cones $\sigma_I$ with $I \subsetneq
  \{0,1,\dots,n\}$ and $|I| \leq k$, whose maximal cones
  all have weight one (cf. \cite[example 3.9]{AR07}).
  This fan $L^n_k$ is a representative of
  the tropical linear space $\max \{ 0, x_1,\ldots,x_n \}^{n-k}
  \cdot \RR^n$.
\end{definition}

\begin{definition}
  Let $C \in Z_k(\RR^n)$ be a tropical cycle and let the map $i: \RR^n
  \rightarrow \RR^n \times \RR^n$ be given by $x \mapsto (x,x)$.
  Then the push-forward cycle $$ \triangle_C := i_{*}(C) \in
  Z_k(\RR^n \times \RR^n)$$ is called the \emph{diagonal} of $C \times C$.
\end{definition}

In order to express the diagonal in $L^n_k \times L^n_k$ by means
of Cartier divisors we first have to refine $L^n_k \times L^n_k$
in such a way that the diagonal is a subfan of this refinement:

\begin{definition}
  Let $F^n_k$ be the refinement of $L^n_k \times L^n_k$ that arises recursively from
  $L^n_k \times L^n_k$ as follows:
  Let $M := (L^n_k \times L^n_k)^{(2k)}$ be the set of maximal
  cones in $L^n_k \times L^n_k$.
  If a cone $\sigma \in M$ is generated by
  $$\left( \begin{array}{c} -e_{i} \\ \hline 0 \end{array} \right)
  , \left( \begin{array}{c} 0 \\ \hline -e_i \end{array} \right), v_3, \ldots, v_{2k}$$
  for some $i$ and vectors
  $$v_j \in \left\{ \left. \left( \begin{array}{c} -e_\mu \\ \hline -e_\mu \end{array} \right),
  \left( \begin{array}{c} -e_\mu \\ \hline 0 \end{array} \right),
  \left( \begin{array}{c} 0 \\ \hline -e_\mu \end{array} \right) \right| \mu =0,\ldots,n \right\}$$
  then replace the cone $\sigma$ by the two cones spanned by
  $$\left( \begin{array}{c} -e_{i} \\ \hline -e_{i} \end{array} \right),
  \left( \begin{array}{c} -e_{i} \\ \hline 0 \end{array} \right)
  ,v_3,\ldots,v_{2k}$$
  and
  $$\left( \begin{array}{c} -e_{i} \\ \hline -e_{i} \end{array} \right),
  \left( \begin{array}{c} 0 \\ \hline -e_{i} \end{array} \right)
  ,v_3,\ldots,v_{2k},$$
  respectively. Repeat this process until there are no more
  cones in $M$ that can be replaced. The fan $F^n_k$ is then the set of
  all faces of all cones in $M$.
\end{definition}

The next lemma provides a technical tool needed in the proofs of
the subsequent theorems:

\begin{lemma} \label{lem-zerofunctionzeroweight}
  Let $F$ be a complete and smooth fan in $\RR^n$ (in the sense of toric geometry)
  and let the weight of every maximal cone in $F$ be one. Moreover, let $h_1,\ldots,h_r$, $r \leq n$, be
  rational functions on $\RR^n$ that are linear on every cone of
  $F$. Then the intersection product $h_1 \cdots h_r \cdot F$ is
  given by
  $$h_1 \cdots h_r \cdot F = \left( \bigcup_{i=0}^{n-r} F^{(i)},\omega_{h_1 \cdots h_r} \right)$$
  with some weight function $\omega_{h_1 \cdots h_r}$ on the cones
  of dimension $n-r$.\\
  Let $\tau \in F^{(n-r)}$ be a cone in $F$ such that for all
  maximal cones $\sigma \in F^{(n)}$ with $\tau \subseteq \sigma$ there exists
  some index $i \in \{1,\ldots,r \}$ such that $h_i$ is identically zero on $\sigma$. Then
  holds: $$\omega_{h_1 \cdots h_r}(\tau)=0.$$
\end{lemma}
\begin{proof}
  We proof the claim by induction on $r$: For $r=1$ we are in the
  situation that $h_1$ is identically zero on every maximal cone adjacent to
  $\tau$. Hence $\omega_{h_1}(\tau)=0$. Now let $r>1$. Using the
  induction hypothesis we can conclude that $|h_1 \cdots h_{r-1} \cdot
  F| \subseteq \bigcup_{\sigma \in S} \sigma$, where $$S:= {\{\sigma \in F^{(n)} | \text{none of } h_1, \ldots, h_{r-1} \text{ is identically zero on
  } \sigma \}}.$$ Our above assumption then implies that $h_r$ must be identically
  zero on every cone in $$\{ \sigma \in F^{(n)} | \tau \subseteq \sigma \text{ and none of } h_1, \ldots, h_{r-1} \text{ is identically zero on } \sigma
  \}$$ and thus that $\omega_{h_1 \cdots h_r}(\tau)=0$.
\end{proof}

\begin{notation} \label{notation-vector_as_rat_function}
  Let $F$ be a simplicial fan in $\RR^n$ and let $u$ be
  a generator of a ray $r_u$ in $F$. By abuse of notation we also denote by
  $u$ the unique rational function on $|F|$ that is linear on
  every cone in $F$, that has the value one on $u$ and that is
  identically zero on all rays of $F$ other than $r_u$.

  If not stated otherwise, vectors considered as Cartier divisors will from
  now on always denote rational functions on the complete fan $F^n_n$.
\end{notation}

\begin{notation}
  Let $C$ be a tropical cycle and let $h_1,\ldots,h_r \in \Div(C)$
  be Cartier divisors on $C$. If $$P(x_1,\ldots,x_r)=\sum_{i_1+\ldots+i_r \leq d}
  \alpha_{i_1,\ldots,i_r} x_1^{i_1} \cdots x_r^{i_r}$$ is a polynomial
  in variables $x_1,\ldots,x_r$ we denote by $P(h_1,\ldots,h_r) \cdot C$
  the intersection product
  $$P(h_1,\ldots,h_r) \cdot C := \sum_{i_1+\ldots+i_r \leq d}
  \left( \alpha_{i_1,\ldots,i_r} h_1^{i_1} \cdots h_r^{i_r} \cdot C \right).$$
\end{notation}

In the following theorem we give a description of the diagonal
$\triangle_{L^n_{n-k}}$ by means of Cartier divisors on our fan
$F^n_n$:

\begin{theorem}\label{thm-diagonalinRn}
  The fan
  $$ \left( \left( \begin{array}{c} -e_1 \\ \hline 0 \end{array} \right) +
  \left( \begin{array}{c} 0 \\ \hline -e_0 \end{array} \right) \right)
  \dots \left( \left( \begin{array}{c} -e_n \\ \hline 0 \end{array} \right) +
  \left( \begin{array}{c} 0 \\ \hline -e_0 \end{array} \right) \right) \cdot
  \left( \left( \begin{array}{c} -e_0 \\ \hline 0 \end{array} \right) +
  \left( \begin{array}{c} -e_0 \\ \hline -e_0 \end{array} \right) \right)^k
  \cdot F^n_n $$
  is a representative of the diagonal $\triangle_{L^n_{n-k}}$.
\end{theorem}
\begin{proof}
  First of all, note that $$\left( \begin{array}{c} -e_0 \\ \hline 0 \end{array} \right) +
  \left( \begin{array}{c} -e_0 \\ \hline -e_0 \end{array} \right) $$ is
  a representation of the tropical polynomial $\max \{0,x_1,\ldots,x_n \}$, where $x_1,\ldots,x_n$
  are the coordinates of the first factor of $\RR^n \times \RR^n$.
  Applying \cite[lemma 9.6]{AR07} we obtain
  $$\left[ \left( \left( \begin{array}{c} -e_0 \\ \hline 0 \end{array} \right) +
  \left( \begin{array}{c} -e_0 \\ \hline -e_0 \end{array} \right) \right)^k \cdot F^n_n \right]= [L^n_{n-k} \times \RR^n].$$
  By lemma \cite[lemma 9.4]{AR07} we can conclude that
  $\triangle_{\RR^n} \cdot [L^n_{n-k} \times \RR^n] =
  i_{*}([L^n_{n-k}])=\triangle_{L^n_{n-k}}$ and hence it suffices to show
  that $[X]=\triangle_{\RR^n}$ for
  $$ X:= \left( \left( \begin{array}{c} -e_1 \\ \hline 0 \end{array} \right) +
  \left( \begin{array}{c} 0 \\ \hline -e_0 \end{array} \right) \right)
  \dots \left( \left( \begin{array}{c} -e_n \\ \hline 0 \end{array} \right) +
  \left( \begin{array}{c} 0 \\ \hline -e_0 \end{array} \right)
  \right) \cdot F^n_n$$
  to prove the claim. Therefore, let $\sigma = \langle r_1, \ldots, r_n \rangle_{\RR_{\geq 0}} \in X^{(n)}$
  be a cone not contained in $|\triangle_{\RR^n}|$. We will show that the weight of $\sigma$ in $X$ has to be zero.
  W.l.o.g. we assume that $$r_1 \not\in
  D:=\left\{ \left( \begin{array}{c} -e_0 \\ \hline -e_0 \end{array} \right),\dots,
  \left( \begin{array}{c} -e_n \\ \hline -e_n \end{array}
  \right) \right\}.$$
  Moreover, let $$T:= \left\{ \left( \begin{array}{c} -e_1 \\ \hline 0 \end{array} \right),\dots,
  \left( \begin{array}{c} -e_n \\ \hline 0 \end{array}
  \right) \right\} \text{ and } B:=\left\{ \left( \begin{array}{c} 0 \\ \hline -e_1 \end{array} \right),\dots,
  \left( \begin{array}{c} 0 \\ \hline -e_n \end{array}
  \right) \right\}.$$
  We distinguish between two cases:
  \begin{enumerate}
  \item[1.] First, we assume that $$r_i \not\in \left\{\left( \begin{array}{c} -e_0 \\ \hline 0 \end{array} \right),
  \left( \begin{array}{c} 0 \\ \hline -e_0 \end{array} \right)
  \right\}, i=1,\ldots,n. $$
  Changing the given rational functions by globally linear functions we can rewrite the above intersection product as $X=\varphi_1
  \cdots \varphi_n \cdot F^n_n$, where
  $$\varphi_i = \left\{ \begin{array}{ll}
    \left( \begin{array}{c} -e_i \\ \hline 0 \end{array} \right)+\left( \begin{array}{c} 0 \\ \hline -e_0 \end{array}
    \right), & \text{if } \left( \begin{array}{c} -e_i \\ \hline 0 \end{array} \right) \not\in \{r_1,\ldots,r_n\} \vspace{2mm} \\
    \left( \begin{array}{c} 0 \\ \hline -e_i \end{array} \right)+\left( \begin{array}{c} -e_0 \\ \hline 0 \end{array}
    \right), & \text{else.}
  \end{array} \right.$$
  Now we apply lemma \ref{lem-zerofunctionzeroweight}: If the weight of $\sigma$ in $X$
  is non-zero there must be at least one cone $$\widetilde{\sigma} =\langle r_1,\ldots,r_n,v_1,\dots,v_n
  \rangle_{\RR \geq 0} \in F^n_n$$ such that all rational functions $\varphi_1,\ldots,\varphi_n$
  are non-zero on $\widetilde{\sigma}$. We study three subcases:
    \begin{enumerate}
      \item There are vectors $r_i \in
      T$ and $r_j \in B$: Then we need both vectors $\left( \begin{array}{c} -e_0 \\ \hline 0 \end{array}
      \right)$ and $\left( \begin{array}{c} 0 \\ \hline -e_0 \end{array}
      \right)$ among the $v_\mu$ such that all functions $\varphi_i$
      are non-zero on $\widetilde{\sigma}$. But there is no cone in
      $F^n_n$ containing these two vectors.
      \item $r_1 \in T$ (or $r_1 \in B$) and $r_j \in D$ for some $j$ and $r_i \in T \cup D$ (or $r_i \in B \cup D$) for all $i$:
      As there is no cone in $F$ containing $\left( \begin{array}{c} -e_i \\ \hline 0 \end{array}
      \right)$ and $\left( \begin{array}{c} 0 \\ \hline -e_i \end{array}
      \right)$ for any $i$, we need $\left( \begin{array}{c} -e_0 \\ \hline 0 \end{array}
      \right)$ among the $v_\mu$ such that all functions $\varphi_i$ are non-zero on $\widetilde{\sigma}$.
      Moreover, if $\left( \begin{array}{c} -e_i \\ \hline 0 \end{array} \right) \not\in \{r_1,\ldots,r_n\}$ then
      we must have $\left( \begin{array}{c} -e_i \\ \hline 0 \end{array}
      \right) \in \{v_1,\ldots,v_n\}$. But there is no cone in $F^n_n$
      containing $\left( \begin{array}{c} -e_1 \\ \hline 0 \end{array}
      \right), \ldots, \left( \begin{array}{c} -e_n \\ \hline 0 \end{array}
      \right)$ and $\left( \begin{array}{c} -e_0 \\ \hline 0 \end{array}
      \right)$. (Analogously for $B$, but with $\varphi_i$ defined the other way around.)
      \item All vectors $r_i$ are contained in $T$ (or in $B$):
      In this case we need $\left( \begin{array}{c} 0 \\ \hline -e_1 \end{array}
      \right)$ or $\left( \begin{array}{c} -e_0 \\ \hline 0 \end{array}
      \right)$ among the $v_\mu$ such that all functions $\varphi_i$
      are non-zero, but again there is no such cone.
      (Analogously for $B$, but with $\varphi_i$ defined the other way around.)
    \end{enumerate}
    \item[2.] Now we assume that $$r_1= \left( \begin{array}{c} -e_0 \\ \hline 0 \end{array}
    \right) \text{ } \left( \text{or } r_1= \left( \begin{array}{c} 0 \\ \hline -e_0 \end{array}
    \right) \right).$$
    Like before we rewrite the intersection product as $X=\varphi_1
    \cdots \varphi_n \cdot F^n_n$ with $\varphi_i$ defined as
    above and apply lemma \ref{lem-zerofunctionzeroweight}:
    If $\left( \begin{array}{c} -e_i \\ \hline 0 \end{array} \right) \not\in
    \{r_1,\ldots,r_n\}$ then $\varphi_i=\left( \begin{array}{c} -e_i \\ \hline 0 \end{array}
    \right)+\left( \begin{array}{c} 0 \\ \hline -e_0 \end{array}
    \right)$ and we need $\left( \begin{array}{c} -e_i \\ \hline 0 \end{array}
    \right)$ or $\left( \begin{array}{c} 0 \\ \hline -e_0 \end{array}
    \right)$ among the $v_\mu$ such that all functions $\varphi_i$ are non-zero on $\widetilde{\sigma}$.
    But as there is no cone in $F^n_n$ containing
    $\left( \begin{array}{c} 0 \\ \hline -e_0 \end{array} \right)$
    and $\left( \begin{array}{c} -e_0 \\ \hline 0 \end{array}
    \right)$ we must have $\left( \begin{array}{c} -e_i \\ \hline 0 \end{array}
    \right) \in \{v_1,\ldots,v_n\}$. Hence all the vectors $\left( \begin{array}{c} -e_1 \\ \hline 0 \end{array}
    \right),\ldots,\left( \begin{array}{c} -e_n \\ \hline 0 \end{array}
    \right)$ and $\left( \begin{array}{c} -e_0 \\ \hline 0 \end{array}
    \right)$ must be contained in
    $\{r_1,\ldots,r_n,v_1,\ldots,v_n\}$, but there is no such cone
    in $F^n_n$. (Analogously for $r_1= \left( \begin{array}{c} 0 \\ \hline -e_0 \end{array}
    \right)$, but with $\varphi_i$ defined the other way around.)
  \end{enumerate}
  So far we have proven that our intersection cycle $X$ is
  contained in the diagonal $\triangle_{\RR^n}$. As the diagonal is
  irreducible we can then conclude by \cite[lemma 2.21]{GKM07} that
  $[X]=\lambda \cdot \triangle_{\RR^n}$ for some integer $\lambda$.
  Thus our last step in this proof is to show that $\lambda=1$:
  Let $\varphi_1,\ldots,\varphi_n$ be the rational functions given
  above. We obtain the following equation
  of cycles in $\RR^n \times \RR^n$:
  $$\begin{array}{rcl}
    && \varphi_1 \cdots \varphi_n \cdot [\{0\} \times \RR^n] \vspace{2mm}\\
    &=& \left( \left( \begin{array}{c} -e_1 \\ \hline 0 \end{array} \right) +
    \left( \begin{array}{c} 0 \\ \hline -e_0 \end{array} \right) \right)
    \dots \left( \left( \begin{array}{c} -e_n \\ \hline 0 \end{array} \right) +
    \left( \begin{array}{c} 0 \\ \hline -e_0 \end{array} \right)
    \right) \cdot [\{0\} \times \RR^n] \vspace{2mm}\\
    &=& \left( \begin{array}{c} 0 \\ \hline -e_0 \end{array} \right)^n \cdot [\{0\} \times \RR^n] \vspace{2mm}\\
    &=& \{0\} \times \{0\}.
  \end{array}$$
  As $\varphi_1 \cdots \varphi_n \cdot [\RR^n \times
  \RR^n]= \lambda \cdot \triangle_{\RR^n}$, by \cite[definition 9.3]{AR07} and \cite[remark 9.9]{AR07}
  we obtain the equation
  $$\begin{array}{rcl}
    \lambda \cdot  \{0\} &=& \lambda \cdot ( \{0\} \cdot \RR^n)\\
    &=& \pi_{*}(\varphi_1 \cdots \varphi_n \cdot (\{0\} \times \RR^n) )\\
    &=& \pi_{*}(\{0\} \times \{0\}) )\\
    &=& 1 \cdot \{0\}
  \end{array}$$
  of cycles in $\RR^n$. This finishes the proof.
\end{proof}

Our next step is to derive a description of the diagonal
$\triangle_{L^n_{n-k}}$ on $L^n_{n-k} \times L^n_{n-k}$ from our
description on $F^n_n$:

\begin{theorem}\label{thm-diagonalinLnk}
  The intersection product in theorem \ref{thm-diagonalinRn} can
  be rewritten as
  $$ \left( \sum_{i=1}^r h_{i,1} \dots h_{i,n-k} \right) \cdot
  \left( \left( \begin{array}{c} 0 \\ \hline -e_0 \end{array} \right) +
  \left( \begin{array}{c} -e_0 \\ \hline -e_0 \end{array} \right) \right)^k \cdot
  \left( \left( \begin{array}{c} -e_0 \\ \hline 0 \end{array} \right) +
  \left( \begin{array}{c} -e_0 \\ \hline -e_0 \end{array} \right) \right)^k
  \cdot F^n_n$$
  for some Cartier divisors $h_{i,j}$ on $F^n_n$.
\end{theorem}

We have to prepare the proof of the theorem by the following
lemma:

\pagebreak

\begin{lemma} \label{lemma-neededrelations}
  Let $C \in Z_l(L^n_{n-k})$ be a subcycle of $L^n_{n-k}$. Then the
  following intersection products are zero:
  \begin{enumerate}
    \item $\left( \begin{array}{c} -e_0 \\ \hline 0 \end{array} \right) \cdot \left( \begin{array}{c} 0 \\ \hline -e_0 \end{array}
    \right) \cdot (C \times \RR^n)$ \vspace{2mm},
    \item $v_{i_1} \cdots v_{i_{n-k+r}} \cdot (C \times \RR^n)$ \vspace{2mm},
    \item $\left( \begin{array}{c} 0 \\ \hline -e_0 \end{array} \right) \cdot \left( \begin{array}{c} -e_0 \\ \hline -e_0 \end{array}
    \right)^s \cdot v_{i_1} \cdots v_{i_{n-k-s+r}} \cdot (C \times \RR^n),$
  \end{enumerate}
  where $r,s>0$ and the vectors $$v_{i_j} \in \left\{ \left( \begin{array}{c} -e_1 \\ \hline 0 \end{array}
    \right),\ldots,\left( \begin{array}{c} -e_n \\ \hline 0 \end{array} \right),
    \left( \begin{array}{c} -e_0 \\ \hline -e_0 \end{array} \right)\right\}$$ are pairwise distinct.
\end{lemma}
\begin{proof}
  (a) and (b): In both cases, a cone that can occur in the
  intersection product with non-zero weight has to be contained
  in a cone of $F^n_n$ that is contained in $|L^n_{n-k} \times \RR^n|$
  and that contains the vectors $\left( \begin{array}{c}
  -e_0 \\ \hline 0 \end{array} \right), \left( \begin{array}{c} 0
  \\ \hline -e_0 \end{array} \right)$ or $v_{i_1}, \ldots,
  v_{i_{n-k+r}}$, respectively. But there are no such cones.\\
  (c): By (a) and \cite[lemma 9.7]{AR07} we can rewrite the intersection product as
  $$\begin{array}{rl}
    & \left( \begin{array}{c} 0 \\ \hline -e_0 \end{array} \right) \cdot \left( \begin{array}{c} -e_0 \\ \hline -e_0 \end{array}
    \right)^s \cdot v_{i_1} \cdots v_{i_{n-k-s+r}} \cdot (C \times \RR^n)\\
    =& \left( \begin{array}{c} 0 \\ \hline -e_0 \end{array} \right) \cdot \left( \left( \begin{array}{c} -e_0 \\ \hline 0 \end{array}
    \right)+\left( \begin{array}{c} -e_0 \\ \hline -e_0 \end{array}
    \right) \right)^s \cdot v_{i_1} \cdots v_{i_{n-k-s+r}} \cdot (C \times \RR^n)\\
    =& \left( \begin{array}{c} 0 \\ \hline -e_0 \end{array} \right) \cdot v_{i_1} \cdots v_{i_{n-k-s+r}} \cdot
    \left[ \left( \left( \begin{array}{c} -e_0 \\ \hline 0 \end{array}
    \right)+\left( \begin{array}{c} -e_0 \\ \hline -e_0 \end{array}
    \right) \right)^s \cdot C \right] \times \RR^n\\
    =& \left( \begin{array}{c} 0 \\ \hline -e_0 \end{array} \right) \cdot v_{i_1} \cdots v_{i_{n-k-s+r}} \cdot
    \left[ \max \{0,x_1,\ldots,x_n \}^s \cdot C \right] \times \RR^n,
  \end{array}$$
  which is zero by (b) as $\max \{0,x_1,\ldots,x_n \}^s \cdot C$ is contained in $L^n_{n-k-s}$.
\end{proof}

\begin{proof}[Proof of theorem \ref{thm-diagonalinLnk}]
  By theorem \ref{thm-diagonalinRn} we have the representation
  $$\begin{array}{rcl}
    \triangle_{L^n_{n-k}} &=& { \tiny \left( \left( \begin{array}{c} -e_1 \\ \hline 0 \end{array} \right) +
    \left( \begin{array}{c} 0 \\ \hline -e_0 \end{array} \right) \right)
    \dots \left( \left( \begin{array}{c} -e_n \\ \hline 0 \end{array} \right) +
    \left( \begin{array}{c} 0 \\ \hline -e_0 \end{array} \right) \right)} \cdot
    \underbrace{ {\tiny \left( \left( \begin{array}{c} -e_0 \\ \hline 0 \end{array} \right) +
    \left( \begin{array}{c} -e_0 \\ \hline -e_0 \end{array} \right)
    \right)^k} \cdot [F^n_n]}_{=[L^n_{n-k} \times \RR^n]} \\
    &=& { \tiny \left( \left( \begin{array}{c} -e_1 \\ \hline 0 \end{array} \right)
    \cdots
    \left( \begin{array}{c} -e_n \\ \hline 0 \end{array} \right)
    + \ldots + \left( \begin{array}{c} 0 \\ \hline -e_0 \end{array} \right)^n \right) } \cdot
    [L^n_{n-k} \times \RR^n].
  \end{array}$$
  By lemma \ref{lemma-neededrelations} (b) all the summands containing $\left( \begin{array}{c} 0 \\ \hline -e_0 \end{array}
  \right)^s$ with a power $s<k$ are zero. Hence we can rewrite the
  intersection product as
  {\small $$\begin{array}{rcl} \triangle_{L^n_{n-k}} &=&  \left[ \left( \begin{array}{c} -e_1 \\ \hline 0 \end{array} \right)\cdots
    \left( \begin{array}{c} -e_{n-k} \\ \hline 0 \end{array} \right)
    + \ldots + \left( \begin{array}{c} 0 \\ \hline -e_0 \end{array} \right)^{n-k} \cdot
    \left( \left( \begin{array}{c} 0 \\ \hline -e_0 \end{array} \right) + \left( \begin{array}{c} -e_0 \\ \hline -e_0 \end{array} \right)
    \right)^k  -A \right] \vspace{1mm}\\ & & \cdot [ L^n_{n-k} \times \RR^n],
  \end{array}$$}where $A$ contains all the summands we added too much.
  Thus all the summands of $A$ are of the form
  $$ \alpha \cdot v_1 \cdots v_{n-s-t} \cdot \left( \begin{array}{c} 0 \\ \hline -e_0 \end{array}
  \right)^{s} \cdot \left( \begin{array}{c} -e_0 \\ \hline -e_0 \end{array}
  \right)^t$$
  for some integer $\alpha$, vectors $v_{i} \in \left\{ \left( \begin{array}{c} -e_1 \\ \hline 0 \end{array} \right),\dots,
  \left( \begin{array}{c} -e_n \\ \hline 0 \end{array}
  \right) \right\}$ and powers $1 \leq t \leq k$, ${0 \leq s \leq n}$. By lemma
  \ref{lemma-neededrelations} (b) and (c) such a summand applied to $[ L^n_{n-k} \times \RR^n]$ is zero if $s<k$
  and only those summands remain in $A$ that have $t \geq 1, s \geq
  k$. Let $$ S:=\alpha \cdot v_1 \cdots v_{n-s-t} \cdot \left( \begin{array}{c} 0 \\ \hline -e_0 \end{array}
  \right)^{s} \cdot \left( \begin{array}{c} -e_0 \\ \hline -e_0 \end{array}
  \right)^t$$ be one of the remaining summands. By lemma
  \ref{lemma-neededrelations} (a) we obtain the equation
  $$\begin{array}{rl}
    & \alpha \cdot v_1 \cdots v_{n-s-t} \cdot { \tiny \left( \begin{array}{c} 0 \\ \hline -e_0 \end{array}
    \right)^{s} \cdot \left( \begin{array}{c} -e_0 \\ \hline -e_0 \end{array}
    \right)^t \cdot [ L^n_{n-k} \times \RR^n] }\\
    = & {\tiny \left( \sum\limits_{j=0}^t {{t}\choose{j}} \cdot \alpha \cdot v_1 \cdots v_{n-s-t} \cdot \left( \begin{array}{c} 0 \\ \hline -e_0 \end{array}
    \right)^{s} \cdot \left( \begin{array}{c} -e_0 \\ \hline -e_0 \end{array}
    \right)^j \cdot \left( \begin{array}{c} -e_0 \\ \hline 0 \end{array}
    \right)^{t-j} \right) } \cdot [ L^n_{n-k} \times \RR^n]\\
    = & {\tiny \left( \alpha \cdot v_1 \cdots v_{n-s-t} \cdot \left( \begin{array}{c} 0 \\ \hline -e_0 \end{array}
    \right)^{s} \cdot \left( \left( \begin{array}{c} -e_0 \\ \hline 0 \end{array}
    \right) + \left( \begin{array}{c} -e_0 \\ \hline -e_0 \end{array}
    \right) \right)^t \right) } \cdot [ L^n_{n-k} \times \RR^n]\\
    = & { \tiny \left[ \left( \left( \begin{array}{c} 0 \\ \hline -e_0 \end{array} \right) + \left( \begin{array}{c} -e_0 \\ \hline -e_0 \end{array} \right) \right)^k
    \right. } \\ & { \tiny \left. \cdot \left( \alpha \cdot v_1 \cdots v_{n-s-t} \cdot \left( \begin{array}{c} 0 \\ \hline -e_0 \end{array}
    \right)^{s-k} \cdot \left( \left( \begin{array}{c} -e_0 \\ \hline 0 \end{array}
    \right) + \left( \begin{array}{c} -e_0 \\ \hline -e_0 \end{array}
    \right) \right)^t \right) -B_S \right] } \cdot [ L^n_{n-k} \times \RR^n],
  \end{array}$$
  where $B_S$ contains again all the summands we added too much.
  Thus all the summands of $B_S$ are of the form
  $$ S':=\beta \cdot {{t}\choose{t'}} \cdot v_1 \cdots v_{n-s-t} \cdot \left( \begin{array}{c} 0 \\ \hline -e_0 \end{array}
    \right)^{s-s'} \cdot \left( \begin{array}{c} -e_0 \\ \hline -e_0 \end{array}
    \right)^{s'} \cdot \left( \begin{array}{c} -e_0 \\ \hline 0 \end{array}
    \right)^{t'} \cdot \left( \begin{array}{c} -e_0 \\ \hline -e_0 \end{array}
    \right)^{t-t'}$$
  for some integer $\beta$ and powers $1 \leq s' \leq k$, $0 \leq t' \leq t$. If $s-s'<k$ we group all corresponding summands together as
  $$ \beta \cdot v_1 \cdots v_{n-s-t} \cdot \left( \begin{array}{c} 0 \\ \hline -e_0 \end{array}
  \right)^{s-s'} \cdot \left( \begin{array}{c} -e_0 \\ \hline -e_0 \end{array}
  \right)^{s'} \cdot \left( \left( \begin{array}{c} -e_0 \\ \hline 0 \end{array}
  \right) + \left( \begin{array}{c} -e_0 \\ \hline -e_0 \end{array}
  \right) \right)^t. $$
  This product applied to $[L^n_{n-k} \times \RR^n]$ is zero by lemma \ref{lemma-neededrelations} (b) and
  (c). Moreover, all summands $S'$ with $s-s' \geq k$ and $t'>0$
  yield zero on $[ L^n_{n-k} \times \RR^n]$ by lemma \ref{lemma-neededrelations} (a). Thus only those summands $S'$ are left in $B_S$ that are of the form
  $$ S'=\beta' \cdot v_1 \cdots v_{n-s-t} \cdot \left( \begin{array}{c} 0 \\ \hline -e_0 \end{array}
    \right)^{s-s'} \cdot \left( \begin{array}{c} -e_0 \\ \hline -e_0 \end{array}
    \right)^{t+s'}$$
  with $s-s' \geq k$ and $s' \geq 1$. Applying this process inductively to all summands with $t=1,\ldots,n-k-1$ in which we could not factor out
  $\left( \left( \begin{array}{c} 0 \\ \hline -e_0 \end{array} \right) + \left( \begin{array}{c} -e_0 \\ \hline -e_0 \end{array} \right)
  \right)^k$, yet, we can by and by increase the power of $\left( \begin{array}{c} -e_0 \\ \hline -e_0 \end{array} \right)$
  in all remaining summands until finally only one summand
  $$ \gamma \cdot \left( \begin{array}{c} 0 \\ \hline -e_0 \end{array}
    \right)^{k} \cdot \left( \begin{array}{c} -e_0 \\ \hline -e_0 \end{array}
    \right)^{n-k} $$
  is left. But
  $$\begin{array}{l} \gamma \cdot \left( \begin{array}{c} 0 \\ \hline -e_0 \end{array}
    \right)^{k} \cdot \left( \begin{array}{c} -e_0 \\ \hline -e_0 \end{array}
    \right)^{n-k} \cdot [ L^n_{n-k} \times \RR^n] \\
    = \gamma \cdot \left( \left( \begin{array}{c} 0 \\ \hline -e_0 \end{array}
    \right) + \left( \begin{array}{c} -e_0 \\ \hline -e_0 \end{array}
    \right) \right)^k \cdot \left( \left( \begin{array}{c} -e_0 \\ \hline 0 \end{array}
    \right) + \left( \begin{array}{c} -e_0 \\ \hline -e_0 \end{array}
    \right) \right)^{n-k} \cdot [ L^n_{n-k} \times \RR^n] \end{array}$$
  as
    $$\begin{array}{l} \left( \begin{array}{c} 0 \\ \hline -e_0 \end{array} \right)^{i} \cdot \left( \begin{array}{c} -e_0 \\ \hline -e_0 \end{array}
    \right)^{k-i} \cdot \left( \left( \begin{array}{c} -e_0 \\ \hline 0 \end{array}
    \right) + \left( \begin{array}{c} -e_0 \\ \hline -e_0 \end{array}
    \right) \right)^{n-k} \cdot [L^n_{n-k} \times \RR^n]\\
    =\left( \begin{array}{c} 0 \\ \hline -e_0 \end{array} \right)^{i} \cdot \left( \begin{array}{c} -e_0 \\ \hline -e_0 \end{array}
    \right)^{k-i} \cdot [L^n_{0} \times \RR^n]\\
    = 0 \end{array}$$
  for all $i<k$ by lemma \ref{lemma-neededrelations} (b) and
    $$ \left( \begin{array}{c} 0 \\ \hline -e_0 \end{array} \right)^{k} \cdot \left( \begin{array}{c} -e_0 \\ \hline 0 \end{array} \right)^{j}
    \cdot \left( \begin{array}{c} -e_0 \\ \hline -e_0 \end{array}
    \right)^{n-k-j} \cdot [L^n_{n-k} \times \RR^n] = 0$$
  for all $j>0$ by lemma \ref{lemma-neededrelations} (a). This
  proves the claim.
\end{proof}

\begin{example}
  We perform the steps described in the proof of theorem
  \ref{thm-diagonalinLnk} for the case $n=3, k=2$:\\
  By theorem \ref{thm-diagonalinRn} we have the representation:
$$\begin{array}{rcl}
    \triangle_{L^3_1} &=& { \tiny \left( \left( \begin{array}{c} -e_1 \\ \hline 0 \end{array} \right) +
    \left( \begin{array}{c} 0 \\ \hline -e_0 \end{array} \right) \right)
    \cdot \left( \left( \begin{array}{c} -e_2 \\ \hline 0 \end{array} \right) +
    \left( \begin{array}{c} 0 \\ \hline -e_0 \end{array} \right) \right)
    \cdot \left( \left( \begin{array}{c} -e_3 \\ \hline 0 \end{array} \right) +
    \left( \begin{array}{c} 0 \\ \hline -e_0 \end{array} \right)
    \right)} \vspace{2mm}
    \\ && \cdot \underbrace{ {\tiny \left( \left( \begin{array}{c} -e_0 \\ \hline 0 \end{array} \right) +
    \left( \begin{array}{c} -e_0 \\ \hline -e_0 \end{array} \right)
    \right)^2} \cdot [F^3_3]}_{=[L^3_1 \times \RR^3]} \\
    &=& \Big( { \tiny \underbrace{\left( \begin{array}{c} -e_1 \\ \hline 0 \end{array} \right)
    \cdot \left( \begin{array}{c} -e_2 \\ \hline 0 \end{array} \right)
    \cdot \left( \begin{array}{c} -e_3 \\ \hline 0 \end{array} \right)}_{=0 \text{ by lemma \ref{lemma-neededrelations} (b)}}+
    \underbrace{\left( \begin{array}{c} -e_1 \\ \hline 0 \end{array} \right)
    \cdot \left( \begin{array}{c} -e_2 \\ \hline 0 \end{array} \right)
    \cdot \left( \begin{array}{c} 0 \\ \hline -e_0 \end{array} \right)}_{=0 \text{ by lemma \ref{lemma-neededrelations} (b)}}} \vspace{2mm}\\
    &&+
    { \tiny \underbrace{\left( \begin{array}{c} -e_1 \\ \hline 0 \end{array} \right)
    \cdot \left( \begin{array}{c} -e_3 \\ \hline 0 \end{array} \right)
    \cdot \left( \begin{array}{c} 0 \\ \hline -e_0 \end{array} \right)}_{=0 \text{ by lemma \ref{lemma-neededrelations} (b)}}+
    \underbrace{\left( \begin{array}{c} -e_2 \\ \hline 0 \end{array} \right)
    \cdot \left( \begin{array}{c} -e_3 \\ \hline 0 \end{array} \right)
    \cdot \left( \begin{array}{c} 0 \\ \hline -e_0 \end{array} \right)}_{=0 \text{ by lemma \ref{lemma-neededrelations} (b)}}} \vspace{2mm}\\
    &&+
    { \tiny \left( \begin{array}{c} -e_1 \\ \hline 0 \end{array} \right)
    \cdot \left( \begin{array}{c} 0 \\ \hline -e_0 \end{array} \right)^2+
    \left( \begin{array}{c} -e_2 \\ \hline 0 \end{array} \right)
    \cdot \left( \begin{array}{c} 0 \\ \hline -e_0 \end{array} \right)^2+
    \left( \begin{array}{c} -e_3 \\ \hline 0 \end{array} \right)
    \cdot \left( \begin{array}{c} 0 \\ \hline -e_0 \end{array} \right)^2+
    \left( \begin{array}{c} 0 \\ \hline -e_0 \end{array} \right)^3 } \Big) \vspace{2mm}\\
    && \cdot [L^3_1 \times \RR^3].
  \end{array}$$
  Now we factor out ${ \tiny \left( \left( \begin{array}{c} 0 \\ \hline -e_0 \end{array} \right)+
    \left( \begin{array}{c} -e_0 \\ \hline -e_0 \end{array} \right)
    \right)^2}$ and subtract all summands we do not need:
  $$\begin{array}{rcl}
    \triangle_{L^3_1} &=&
    \Big( { \tiny \left( \begin{array}{c} -e_1 \\ \hline 0 \end{array} \right)+
    \left( \begin{array}{c} -e_2 \\ \hline 0 \end{array} \right)+
    \left( \begin{array}{c} -e_3 \\ \hline 0 \end{array} \right)+
    \left( \begin{array}{c} 0 \\ \hline -e_0 \end{array} \right)
    \bigg) \cdot \bigg( \left( \begin{array}{c} 0 \\ \hline -e_0 \end{array} \right)+
    \left( \begin{array}{c} -e_0 \\ \hline -e_0 \end{array} \right)
    \Big)^2} \cdot [L^3_1 \times \RR^3] \vspace{2mm}\\
    && -\Big( { \tiny \underbrace{\left( \begin{array}{c} -e_1 \\ \hline 0 \end{array} \right)\left( \begin{array}{c} -e_0 \\ \hline -e_0
    \end{array}\right)^2}_{=0 \text{ by \ref{lemma-neededrelations} (b)}}
    +\underbrace{\left( \begin{array}{c} -e_2 \\ \hline 0 \end{array} \right)\left( \begin{array}{c} -e_0 \\ \hline -e_0
    \end{array}\right)^2}_{=0 \text{ by \ref{lemma-neededrelations} (b)}}
    +\underbrace{\left( \begin{array}{c} -e_3 \\ \hline 0 \end{array} \right)\left( \begin{array}{c} -e_0 \\ \hline -e_0 \end{array}\right)^2}_{=0 \text{ by \ref{lemma-neededrelations} (b)}}} \vspace{2mm}\\
    && \tiny{ +\underbrace{\left( \begin{array}{c}  0 \\ \hline -e_0 \end{array} \right)\left( \begin{array}{c} -e_0 \\ \hline -e_0
    \end{array}\right)^2}_{=0 \text{ by \ref{lemma-neededrelations} (c)}}
    +\underbrace{2\left( \begin{array}{c} -e_1 \\ \hline 0 \end{array} \right)\left( \begin{array}{c} 0 \\ \hline -e_0 \end{array}\right)\left( \begin{array}{c} -e_0 \\ \hline -e_0 \end{array}
    \right)}_{=0 \text{ by \ref{lemma-neededrelations} (b)}}
    +\underbrace{2\left( \begin{array}{c} -e_2 \\ \hline 0 \end{array} \right)\left( \begin{array}{c} 0 \\ \hline -e_0 \end{array}\right)\left( \begin{array}{c} -e_0 \\ \hline -e_0 \end{array}
    \right)}_{=0 \text{ by \ref{lemma-neededrelations} (b)}} }\vspace{2mm}\\
    && \tiny{ +\underbrace{2\left( \begin{array}{c} -e_3 \\ \hline 0 \end{array} \right)\left( \begin{array}{c} 0 \\ \hline -e_0 \end{array}\right)\left( \begin{array}{c} -e_0 \\ \hline -e_0 \end{array}
    \right)}_{=0 \text{ by \ref{lemma-neededrelations} (b)}}
    +2\left( \begin{array}{c} 0 \\ \hline -e_0 \end{array} \right)^2\left( \begin{array}{c} -e_0 \\ \hline -e_0 \end{array} \right)
    }\Big) \cdot [L^3_1 \times \RR^3].
  \end{array}$$
  But by lemma \ref{lemma-neededrelations} (a) and (b) we have the following equation for this last summand:
  $$\begin{array}{rl}
    & {\tiny -2\left( \begin{array}{c} 0 \\ \hline -e_0 \end{array} \right)^2\left( \begin{array}{c} -e_0 \\ \hline -e_0 \end{array} \right)}
    \cdot [L^3_1 \times \RR^3] \vspace{2mm}\\
    =&
    {\tiny -2 \left( \left( \begin{array}{c} 0 \\ \hline -e_0 \end{array} \right)^2
    +2\left( \begin{array}{c} 0 \\ \hline -e_0 \end{array} \right)\left( \begin{array}{c} -e_0 \\ \hline -e_0 \end{array} \right)
    +\left( \begin{array}{c} -e_0 \\ \hline -e_0 \end{array} \right)^2
    \right) \cdot \left( \left( \begin{array}{c} -e_0 \\ \hline 0 \end{array} \right) + \left( \begin{array}{c} -e_0 \\ \hline -e_0 \end{array} \right)\right)
    } \cdot [L^3_1 \times \RR^3].
  \end{array}$$
  Hence we obtain altogether:
  $$\begin{array}{rcl}
    \triangle_{L^3_1} &=&
    \Big( { \tiny \left( \begin{array}{c} -e_1 \\ \hline 0 \end{array} \right)+
    \left( \begin{array}{c} -e_2 \\ \hline 0 \end{array} \right)+
    \left( \begin{array}{c} -e_3 \\ \hline 0 \end{array} \right)+
    \left( \begin{array}{c} 0 \\ \hline -e_0 \end{array} \right)-2
    \left( \begin{array}{c} -e_0 \\ \hline 0 \end{array} \right)-2
    \left( \begin{array}{c} -e_0 \\ \hline -e_0 \end{array} \right)
    \bigg)} \vspace{2mm}\\
    && \cdot {\tiny \bigg( \left( \begin{array}{c} 0 \\ \hline -e_0 \end{array} \right)+
    \left( \begin{array}{c} -e_0 \\ \hline -e_0 \end{array} \right)
    \Big)^2} \cdot {\tiny \bigg( \left( \begin{array}{c} -e_0 \\ \hline 0 \end{array} \right)+
    \left( \begin{array}{c} -e_0 \\ \hline -e_0 \end{array} \right)
    \Big)^2} \cdot [\RR^3 \times \RR^3].
  \end{array}$$
\end{example}

\begin{corollary} \label{coro-diagonalinLnk}
  The Cartier divisors $h_{i,j}$ from theorem \ref{thm-diagonalinLnk} provide the following description
  of the diagonal $\triangle_{L^n_{n-k}}$:
  $$\triangle_{L^n_{n-k}} = \sum_{i=1}^r h_{i,1} \dots h_{i,n-k} \cdot
  [L^n_{n-k} \times L^n_{n-k}].$$
\end{corollary}
\begin{proof}
  Let $x_1,\ldots,x_n$ be the coordinates of the first and $y_1,\ldots,y_n$ be coordinates
  of the second factor of $\RR^n \times \RR^n$.
  Applying \cite[lemma 9.6]{AR07} we can conclude that
  $$\begin{array}{l}
  \left[ \left( \left( \begin{array}{c} 0 \\ \hline -e_0 \end{array} \right) +
  \left( \begin{array}{c} -e_0 \\ \hline -e_0 \end{array} \right) \right)^k \cdot
  \left( \left( \begin{array}{c} -e_0 \\ \hline 0 \end{array} \right) +
  \left( \begin{array}{c} -e_0 \\ \hline -e_0 \end{array} \right) \right)^k
  \cdot F^n_n \right] \vspace{2mm} \\ = \left[ \max \{ 0, x_1,\ldots,x_n \}^{k} \cdot \max \{ 0, y_1,\ldots,y_n \}^{k}
  \cdot F^n_n \right] \vspace{2mm} \\ = [L^n_{n-k} \times L^n_{n-k}]
  \end{array}$$
  and hence by theorem \ref{thm-diagonalinRn} and theorem \ref{thm-diagonalinLnk} that $$\sum_{i=1}^r h_{i,1} \dots h_{i,n-k} \cdot
  [L^n_{n-k} \times L^n_{n-k}] = \triangle_{L^n_{n-k}}.$$
\end{proof}

\begin{remark} \label{remark-getdiagonal}
  As lemma \ref{lemma-neededrelations} does not only hold on
  $L^n_{n-k} \times \RR^n$ but also on any $C \times \RR^n$ with $C$ a
  subcycle of $L^n_{n-k}$, the proof of theorem \ref{thm-diagonalinLnk} indeed shows
  that
  $$\begin{array}{rl}
    & \left( \sum\limits_{i=1}^r h_{i,1} \dots h_{i,n-k} \right) \cdot
    \left( \left( \begin{array}{c} 0 \\ \hline -e_0 \end{array} \right) +
    \left( \begin{array}{c} -e_0 \\ \hline -e_0 \end{array} \right) \right)^k
    \cdot (C \times \RR^n) \vspace{1mm}\\
    = & \left( \left( \begin{array}{c} -e_1 \\ \hline 0 \end{array} \right) +
    \left( \begin{array}{c} 0 \\ \hline -e_0 \end{array} \right) \right)
    \dots \left( \left( \begin{array}{c} -e_n \\ \hline 0 \end{array} \right) +
    \left( \begin{array}{c} 0 \\ \hline -e_0 \end{array} \right) \right)
    \cdot (C \times \RR^n)\\
  \end{array}$$
  for all cycles $C \in Z_l(L^n_{n-k})$. Using \cite[corollary 9.8]{AR07} we can conclude that
  $$\begin{array}{rl}
    & \left( \sum\limits_{i=1}^r h_{i,1} \dots h_{i,n-k} \right) \cdot
    \left( \left( \begin{array}{c} 0 \\ \hline -e_0 \end{array} \right) +
    \left( \begin{array}{c} -e_0 \\ \hline -e_0 \end{array} \right) \right)^k
    \cdot (C \times \RR^n) \vspace{1mm}\\
    = & \triangle_{\RR^n} \cdot (C \times \RR^n) \vspace{1mm}\\
    = & \triangle_{C}
  \end{array}$$
  for all such cycles $C$.
\end{remark}

\begin{corollary} \label{coro-diagonalinU}
  Let $\sigma \in L^n_{n-k}$, let $x \in \sigma$ and let $U \subseteq S_{\sigma}=
  \bigcup_{\sigma' \in L^n_{n-k}: \sigma' \supseteq \sigma} (\sigma')^{ri}$
  be an open subset of $|L^n_{n-k}|$ containing $x$. Moreover, let $F$ be the open
  fan $F:=\{ -x+\sigma \cap U| \sigma \in L^n_{n-k} \}$ and
  $\widetilde{F}$ the associated tropical fan. Then there are
  Cartier divisors $h'_{i,j}$ on $\widetilde{F} \times \widetilde{F}$ such that
  $$\triangle_{[\widetilde{F}]} = \sum_{i=1}^r h'_{i,1} \dots h'_{i,n-k} \cdot
  [\widetilde{F} \times \widetilde{F}].$$
\end{corollary}
\begin{proof}
  To obtain the Cartier divisors $h'_{i,j}$ we just have to
  restrict the Cartier divisors $h_{i,j}$ from corollary \ref{coro-diagonalinLnk}
  to the open set $U \times U$, translate them suitably and extend them from $F \times F$ to the
  associated tropical fan $\widetilde{F} \times \widetilde{F}$.
\end{proof}

\pagebreak

\begin{example}
  The following figure shows two fans associated to open subsets
  of $L^3_2$ as in corollary \ref{coro-diagonalinU}:
  \begin{center}
    \includegraphics[scale=0.66]{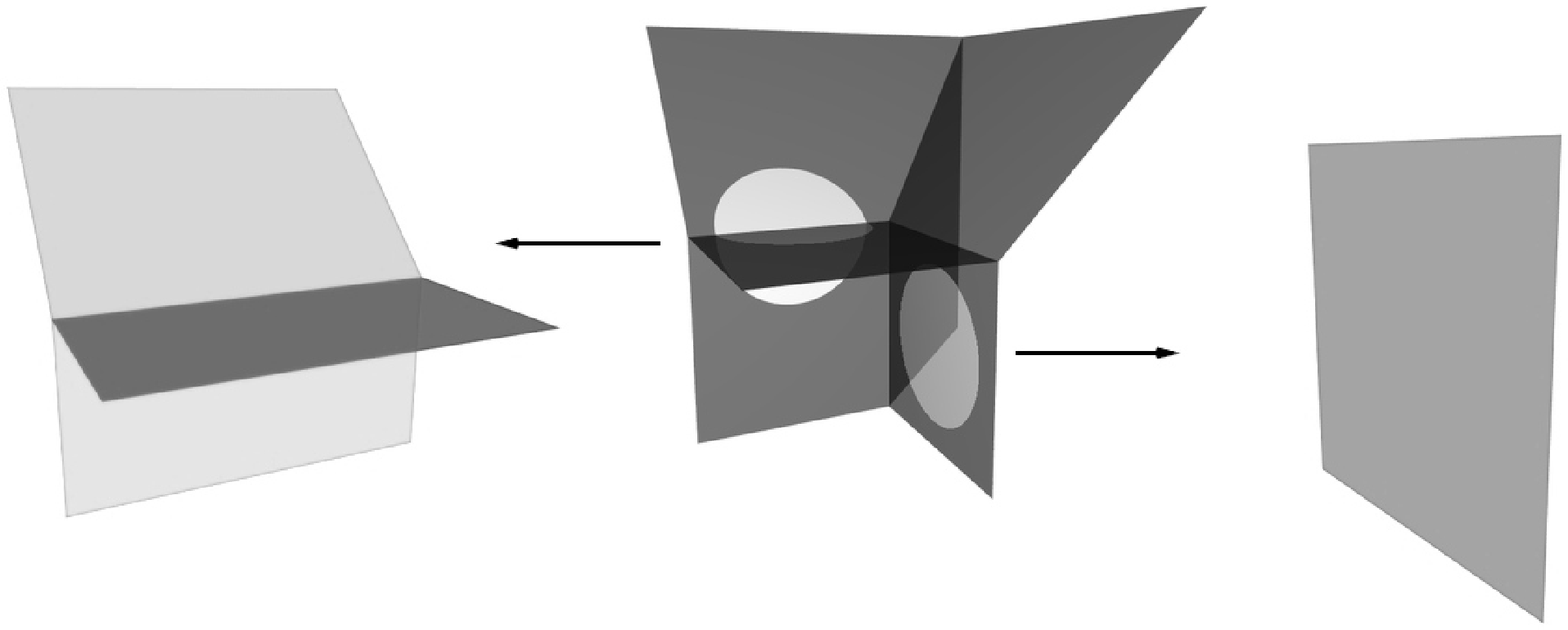}
  \end{center}
\end{example}

\begin{lemma}
  Let $C \in Z_k(\RR^n)$ and $D \in Z_l(\RR^n)$ be tropical cycles such that there exist
  representations of the diagonals $\triangle_C$ and $\triangle_D$ as sums of products of Cartier
  divisors on $C \times C$ and $D \times D$, respectively. Then there also exists a representation of
  $\triangle_{C \times D}$ as a sum of products of Cartier divisors on $(C \times D)^2$.
\end{lemma}
\begin{proof}
  Let $\triangle_C = \sum_{i=1}^r \varphi_{i,1} \dots \varphi_{i,k} \cdot (C \times
  C)$ and $\triangle_D = \sum_{i=1}^s \psi_{i,1} \dots \psi_{i,l} \cdot (D \times
  D)$. \linebreak Moreover, let $\pi_{x}, \pi_{y} :(\RR^n )^4 \rightarrow (\RR^n )^2$
  be given by $(x_1,y_1,x_2,y_2) \mapsto (x_1,x_2)$ and \linebreak $(x_1,y_1,x_2,y_2) \mapsto
  (y_1,y_2)$, respectively.
  Then we have the equation
  $$\triangle_{C \times D} = \left(\sum_{i=1}^r \pi_{x}^{*}\varphi_{i,1} \dots \pi_{x}^{*}\varphi_{i,k} \right)
  \cdot \left(\sum_{i=1}^s \pi_{y}^{*}\psi_{i,1} \dots \pi_{y}^{*}\psi_{i,l} \right)
  \cdot (C \times D)^2.$$
\end{proof}

Now we are ready to define intersection products on all
spaces on which we can express the diagonal by means of Cartier
divisors:

\begin{definition}[Intersection products] \label{def-intersectionproduct}
  Let $C \in Z_k(\RR^n)$ be a tropical cycle and assume that there are Cartier divisors $\varphi_{i,j}$ on
  $C \times C$ such that $$\triangle_{C} = \sum_{i=1}^r \varphi_{i,1} \dots \varphi_{i,k} \cdot
  (C \times C).$$
  Moreover, let $\pi: C \times C \rightarrow C$ be the morphism given by $(x,y) \mapsto x$.
  Then we define the intersection product of subcycles of $C$ by
  $$\begin{array}{rcl}
    Z_{k-l}(C) \times Z_{k-l'}(C) & \longrightarrow & Z_{k-l-l'}(C)\\
    (D_1,D_2) & \longmapsto & D_1 \cdot D_2 := \pi_{*}\left(\sum_{i=1}^r \varphi_{i,1} \dots \varphi_{i,k} \cdot (D_1 \times D_2)
    \right).
  \end{array}$$
\end{definition}

We use the rest of this section to prove that this intersection
product is independent of the used representation of the diagonal
and that it has all the properties we expect --- at least for
those spaces we are interested in:

\begin{lemma}
  Let $C \in Z_k(\RR^n)$ be a tropical cycle,
  $D \in Z_{k-l}(C)$, $E \in Z_{k-l'}(C)$ be subcycles, let $\varphi \in \Div(C)$
  be a Cartier divisor and ${\pi: C \times C \rightarrow C}$ the morphism
  given by $(x,y) \mapsto x$. Then the following equation holds:
  $$ (\varphi \cdot D) \times E = \pi^{*}\varphi \cdot (D \times E).$$
\end{lemma}
\begin{proof}
  The proof is exactly the same as for \cite[lemma 9.6]{AR07}.
\end{proof}

\begin{corollary} \label{coro-cartdivintprodpermute}
  Let $C \in Z_k(\RR^n)$ be a tropical cycle that admits an
  intersection product as in definition
  \ref{def-intersectionproduct}, let $D \in Z_{k-l}(C)$, $E \in
  Z_{k-l'}(C)$ be subcycles and let $\varphi \in \Div(C)$
  be a Cartier divisor. Then the following equation holds:
  $$ (\varphi \cdot D) \cdot E = \varphi \cdot (D \cdot E).$$
\end{corollary}
\begin{proof}
  The proof is exactly the same as for \cite[lemma 9.7]{AR07}.
\end{proof}

\begin{proposition}
  Let $D \in Z_{l}(L^n_{n-k})$ be a subcycle. Then the
  equation $$ D \cdot [L^n_{n-k}] = [L^n_{n-k}] \cdot D = D$$
  holds on $L^n_{n-k}$.
\end{proposition}
\begin{proof}
  Let $\pi_i: L^n_{n-k} \times L^n_{n-k} \rightarrow L^n_{n-k}$ be the morphism given by $(x_1,x_2)
  \mapsto x_i$. By remark \ref{remark-getdiagonal} we get the equation
  $$\begin{array}{rcl}
    D \cdot [L^n_{n-k}] &=& (\pi_1)_{*} \left( \sum\limits_{i=1}^r h_{i,1} \dots h_{i,n-k} \cdot \left( D \times [L^n_{n-k}] \right) \right)\\
    &=& (\pi_1)_{*} \left( \left( \sum\limits_{i=1}^r h_{i,1} \dots h_{i,n-k} \right) \cdot
    \left( \left( \begin{array}{c} 0 \\ \hline -e_0 \end{array} \right) +
    \left( \begin{array}{c} -e_0 \\ \hline -e_0 \end{array} \right) \right)^k
    \cdot (D \times \RR^n) \right)\\
    &=& (\pi_1)_{*} \left( \triangle_{\RR^n} \cdot (D \times \RR^n) \right)\\
    &=& (\pi_1)_{*} \left( \triangle_{D} \right)\\
    &=& D.
  \end{array}$$
  Furthermore, let $\phi:L^n_k \times L^n_k \rightarrow L^n_k \times L^n_k$ be the morphism
  induced by $(x,y) \mapsto (y,x)$. Obviously we have the
  equality
  $$\left( \sum_{i=1}^r h_{i,1} \dots h_{i,n-k} \right) \cdot [L^n_{n-k} \times L^n_{n-k}]
  = \left( \sum_{i=1}^r \phi^{*}h_{i,1} \dots \phi^{*}h_{i,n-k} \right) \cdot [L^n_{n-k} \times L^n_{n-k}].$$
  If $\pi_{ij}: (L^n_{n-k})^4 \rightarrow (L^n_{n-k})^2$
  is the morphism given by $(x_1,x_2,x_3,x_4) \mapsto (x_i,x_j)$ and
  $$ \triangle := \left(\sum_{i=1}^r \pi_{13}^{*}h_{i,1} \dots \pi_{13}^{*}h_{i,n-k} \right)
  \cdot \left(\sum_{i=1}^r \pi_{24}^{*}h_{i,1} \dots \pi_{24}^{*}h_{i,n-k}
  \right)$$
  we can conclude by \cite[proposition 7.7]{AR07} and \cite[lemma 9.6]{AR07} that
  $$ \begin{array}{rl}
    & \left( \sum\limits_{i=1}^r \phi^{*} h_{i,1} \dots \phi^{*} h_{i,n-k} \right) \cdot (D \times [L^n_{n-k}])\\
    =& \left( \sum\limits_{i=1}^r \phi^{*} h_{i,1} \dots \phi^{*} h_{i,n-k} \right) \cdot
    \Big( (D \times [L^n_{n-k}]) \cdot ([L^n_{n-k} \times L^n_{n-k}]) \Big)\\
    =& \left( \sum\limits_{i=1}^r \phi^{*} h_{i,1} \dots \phi^{*} h_{i,n-k} \right) \cdot
    (\pi_{12})_{*} \Big( \triangle \cdot \left( (D \times [L^n_{n-k}]) \times ([L^n_{n-k} \times L^n_{n-k}]) \right) \Big)\\
    =& \left( \sum\limits_{i=1}^r \phi^{*} h_{i,1} \dots \phi^{*} h_{i,n-k} \right) \cdot
    (\pi_{12})_{*} \left( \triangle_{D \times [L^n_{n-k}]} \right)\\
    =& \left( \sum\limits_{i=1}^r \phi^{*} h_{i,1} \dots \phi^{*} h_{i,n-k} \right) \cdot
    (\pi_{34})_{*} \left( \triangle_{D \times [L^n_{n-k}]} \right)\\
    =& (\pi_{34})_{*} \left( \left( \sum\limits_{i=1}^r \pi_{34}^{*}\phi^{*} h_{i,1} \dots \pi_{34}^{*}\phi^{*} h_{i,n-k} \right) \cdot
      \triangle \cdot \left( (D \times [L^n_{n-k}]) \times ([L^n_{n-k} \times L^n_{n-k}]) \right) \right)\\
    =& (\pi_{34})_{*} \left( \triangle \cdot (D \times [L^n_{n-k}]) \times
    \left( \left( \sum\limits_{i=1}^r \phi^{*} h_{i,1} \dots \phi^{*} h_{i,n-k} \right) \cdot [L^n_{n-k} \times L^n_{n-k}] \right) \right)\\
    =& (\pi_{34})_{*} \left( \triangle \cdot (D \times [L^n_{n-k}]) \times
    \left( \left( \sum\limits_{i=1}^r h_{i,1} \dots h_{i,n-k} \right) \cdot [L^n_{n-k} \times L^n_{n-k}] \right) \right)\\
    =& \left( \sum\limits_{i=1}^r h_{i,1} \dots h_{i,n-k} \right) \cdot (D \times [L^n_{n-k}]).
  \end{array}$$
  Hence we can deduce that
  $$\begin{array}{rcl}
    D \cdot [L^n_{n-k}] &=& (\pi_1)_{*} \left( \triangle_{D} \right)\\
    &=& (\pi_2)_{*} \left( \triangle_{D} \right)\\
    &=& (\pi_2)_{*} \left( \left( \sum\limits_{i=1}^r h_{i,1} \dots h_{i,n-k} \right) \cdot (D \times [L^n_{n-k}]) \right)\\
    &=& (\pi_2)_{*} \left( \left( \sum\limits_{i=1}^r \phi^{*} h_{i,1} \dots \phi^{*} h_{i,n-k} \right) \cdot (D \times [L^n_{n-k}]) \right)\\
    &=& (\pi_1)_{*} \left( \left( \sum\limits_{i=1}^r h_{i,1} \dots h_{i,n-k} \right) \cdot ([L^n_{n-k}] \times D) \right)\\
    &=& [L^n_{n-k}] \cdot D.
  \end{array}$$
  This proves the claim.
\end{proof}

\begin{remark}
  We can prove in the same way that $[L^n_{n-k} \times L^m_{m-l}] \cdot D = D$
  holds for all subcycles $D$ of $L^n_{n-k} \times L^m_{m-l}$ and
  even that $[L^{n_1}_{{n_1}-{k_1}} \times \ldots \times  L^{n_r}_{{n_r}-{k_r}}] \cdot D = D$
  holds for all $r \geq 1$ and all subcycles $D$ of $L^{n_1}_{{n_1}-{k_1}} \times \ldots \times  L^{n_r}_{{n_r}-{k_r}}$.
  Moreover, restricting the intersection products to open subsets of $|L^n_k|$ or
  $|L^{n_1}_{{n_1}-{k_1}} \times \ldots \times  L^{n_r}_{{n_r}-{k_r}}|$,
  respectively, implies that $X \cdot D = D$ also holds for all subcycles
  $D \in Z_{l}(X)$ if $X \in \{ [\widetilde{F}], [\widetilde{F_1} \times \ldots \times \widetilde{F_r}] \}$
  where $\widetilde{F}$, $\widetilde{F_i}$ are tropical fans associated to an open
  subsets of some $|L^n_k|$ like in corollary \ref{coro-diagonalinU}.
\end{remark}

\begin{proposition}
  Let $C \in Z_k(\RR^n)$ be a tropical cycle that admits an
  intersection product as in definition \ref{def-intersectionproduct}
  and let $D,D' \in Z_l(C)$, $E \in Z_{l'}(C)$ be subcycles. Then the following equation holds:
  $$(D+D') \cdot E = D \cdot E+ D' \cdot E.$$
\end{proposition}
\begin{proof}
  The proof is exactly the same as for \cite[theorem 9.10 (b)]{AR07}.
\end{proof}

\begin{proposition} \label{propo-independence}
  Let $C \in Z_k(\RR^n)$ be a tropical cycle that admits an
  intersection product as in definition \ref{def-intersectionproduct}
  and let $D \in Z_l(C)$ be a subcycle of $C$. Moreover, let $E \in Z_{l'}(C)$
  be a subcycle such that there are Cartier divisors $\psi_{i,j} \in \Div(C)$
  with $$\sum_{i=1}^r \psi_{i,1} \dots \psi_{i,k-l'} \cdot C = E.$$
  If additionally $C \cdot D = D$ holds then $$\sum_{i=1}^r \psi_{i,1} \dots
  \psi_{i,k-l'} \cdot D = E \cdot D.$$
\end{proposition}
\begin{proof}
  The proof is the same as for \cite[corollary 9.8]{AR07}.
\end{proof}

\begin{remark} \label{remark-independence}
  The meaning of proposition \ref{propo-independence} is the following: If $X \in Z_k(\RR^n)$
  is a tropical cycle such that the diagonal $\triangle_{X}$ can be written
  as a sum of products of Cartier divisors as in definition \ref{def-intersectionproduct}
  and additionally $(X \times X) \cdot Y = Y$ is fulfilled for all subcycles
  $Y$ of $X \times X$ then we can apply proposition \ref{propo-independence} with
  $C:=X \times X$ and $E:=\triangle_X$ to deduce that the definition of the intersection product is
  independent of the choice of the Cartier divisors describing the
  diagonal. In particular we have well-defined intersection products on
  $L^n_k$, $L^{n_1}_{k_1} \times \ldots \times L^{n_r}_{k_r}$, $\widetilde{F}$ and
  $\widetilde{F_1} \times \ldots \times \widetilde{F_r}$
  for all tropical fans $\widetilde{F}$, $\widetilde{F_i}$ associated to an open subset
  of some $|L^n_k|$ like in corollary \ref{coro-diagonalinU}.
\end{remark}

\begin{theorem} \label{theorem-intproductproperties}
  Let $C \in Z_k(\RR^n)$ be a tropical cycle that admits an
  intersection product as in definition \ref{def-intersectionproduct}
  such that additionally $(C \times C) \cdot D = D$ is fulfilled
  for all subcycles $D$ of $C \times C$.
  Moreover, let $E,E' \in Z_l(C)$, $F \in Z_{l'}(C)$ and
  $G \in Z_{l''}(C)$ be subcycles. Then the following equations hold:
  \begin{enumerate}
    \item $E \cdot F = F \cdot E$,
    \item $(E \cdot F) \cdot G = E \cdot (F \cdot G).$
  \end{enumerate}
\end{theorem}
\begin{proof}
  The proof is exactly the same as for \cite[theorem 9.10 (a) and (c)]{AR07}.
\end{proof}

We finish this section with an example showing that even curves
intersecting in the expected dimension can have negative
intersections:

\begin{example}
  Let $C, D \in Z_1(L^3_2)$ be the curves shown in the figure. We
  want to compute the intersection $C \cdot D$. By proposition
  \ref{propo-independence} the easiest way to achieve this is to
  write one of the curves as $\psi \cdot [L^3_2]$ for some Cartier divisor
  $\psi$ on $L^3_2$.
  \begin{center}
    \begin{picture}(180,190)
  \includegraphics[scale=0.75]{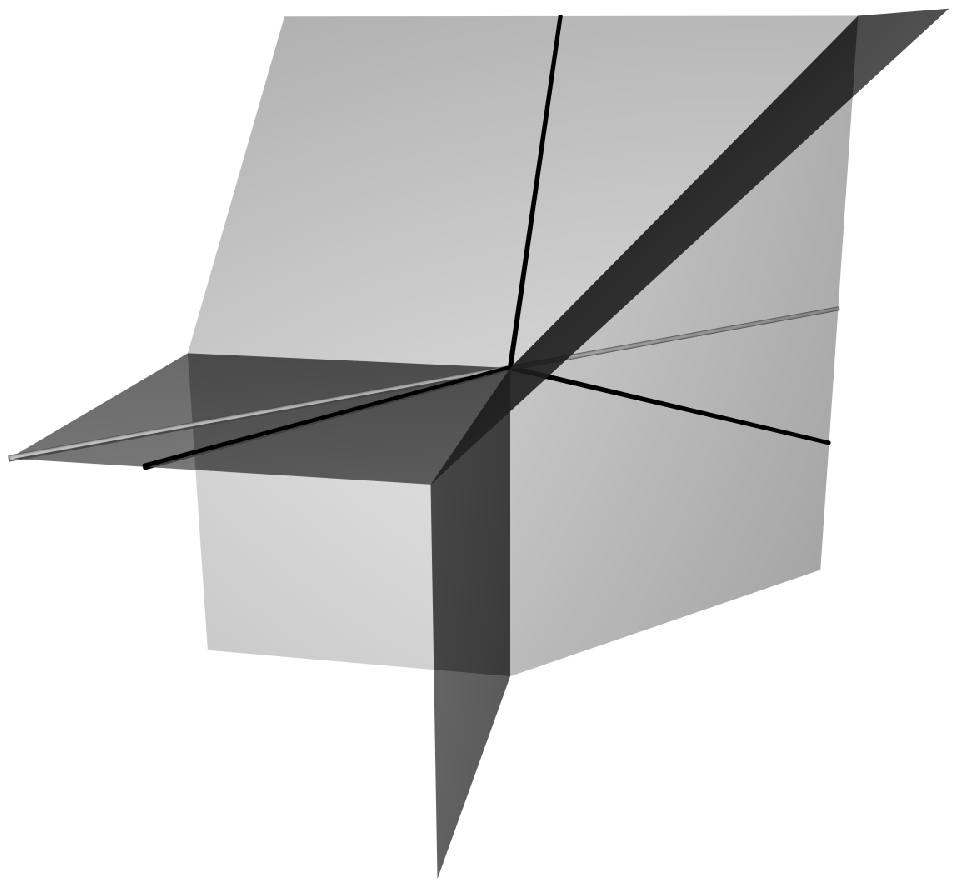}
  \definecolor{curve_C}{gray}{0.3}
  \definecolor{curve_D}{gray}{0.0}
  \put(-205,97){\textcolor{curve_C}{$C$}}
  \put(-237,90){\textcolor{curve_C}{\tiny $\left(\begin{array}{c} -1\\-1\\0\end{array}\right)$}}
  \put(-20,125){\textcolor{curve_C}{\tiny $\left(\begin{array}{c} 1\\1\\0\end{array}\right)$}}
  \put(-86,170){\textcolor{curve_D}{$D$}}
  \put(-195,75){\textcolor{curve_D}{\tiny $\left(\begin{array}{c} -2\\-3\\0\end{array}\right)$}}
  \put(-23,90){\textcolor{curve_D}{\tiny $\left(\begin{array}{c} 2\\2\\-1\end{array}\right)$}}
  \put(-117,170){\textcolor{curve_D}{\tiny $\left(\begin{array}{c} 0\\1\\1\end{array}\right)$}}
\end{picture}

  \end{center}
  Let $F$ be the refinement of $L^3_2$ arising by dividing the cones $\langle -e_1,-e_2 \rangle_{\RR_{\geq 0}}$
  and $\langle -e_0,-e_3 \rangle_{\RR_{\geq 0}}$ into cones $\langle -e_1,-e_1-e_2 \rangle_{\RR_{\geq 0}}$,
  $\langle -e_2,-e_1-e_2 \rangle_{\RR_{\geq 0}}$ and $\langle -e_0,-e_0-e_3 \rangle_{\RR_{\geq 0}}$,
  ${\langle -e_3,-e_0-e_3 \rangle_{\RR_{\geq 0}}}$, respectively. Then
  $$\psi:=\left(\begin{array}{c}1\\1\\1\end{array}\right)-\left(\begin{array}{c}-1\\-1\\0\end{array}\right)$$
  defines a rational function on $F$. As shown in \cite[example 3.10]{AR07}
  we have $\psi \cdot [L^3_2] = C$. Hence we can calculate
  $$\begin{array}{rcl}
    C \cdot D = \psi \cdot D &=& \left( \psi \left(\begin{array}{c}-2\\-3\\0\end{array}\right) + \psi \left(\begin{array}{c}2\\2\\1\end{array}\right)
    + \psi \left(\begin{array}{c}0\\1\\1\end{array}\right) - \psi \left(\begin{array}{c}0\\0\\0\end{array}\right) \right) \cdot \{0\}\\
    &=& (-2+0+1-0) \cdot \{0\}\\
    &=& -1 \cdot \{0\}.
  \end{array}$$
\end{example}

\begin{remark}
  This result is remarkable for the following reason: Our ambient space
  $L^3_2$ arises as a so-called \emph{modification} of $\RR^2$ (cf. \cite{M06}, \cite{M07}).
  Varieties that are connected by a series of modifications are
  called \emph{equivalent} by G. Mikhalkin and are expected to
  have similar properties. But the above example shows that there
  is a big difference between $\RR^2$ and $L^3_2$ even though they
  are equivalent: On $\RR^2$ there is no negative intersection
  product of curves, on $L^3_2$ there is.
\end{remark}

\section{Intersection products on smooth tropical varieties} \label{sec-smoothvarieties}

In this section we use our results from section
\ref{sec-tropicallinearspaces} to define an intersection product
on smooth tropical varieties, i.e. on varieties with tropical
linear spaces as local building blocks:

\begin{definition}[Smooth tropical varieties] \label{def-smoothvarieties}
  An abstract tropical variety $C$ is called a \emph{smooth variety} if it has
  a representative $(((X,|X|),\omega_X),\{ \Phi_\sigma\})$
  such that all the maps
    $$\Phi_\sigma: S_\sigma = \bigcup_{\sigma' \in X^{*}, \sigma' \supset
    \sigma} (\sigma')^{ri} \stackrel{\sim}{\longrightarrow} |F_\sigma| \subseteq |\widetilde{F_\sigma}|$$
  (cf. \cite[definition 5.4]{AR07}) map into tropical fans
  $\widetilde{F_\sigma} = \widetilde{F^\sigma_1} \times \ldots \times \widetilde{F^\sigma_{r_\sigma}}$
  where the $\widetilde{F_i^\sigma}$ are tropical fans associated to open
  subsets of some $|L^{{n_{\sigma,i}}}_{k_{\sigma,i}}|$ as in corollary \ref{coro-diagonalinU}.
\end{definition}

\begin{remark}
  Note that the existence of such a representative $(((X,|X|),\omega_X),\{
  \Phi_\sigma\})$ for $C$ implies that all representatives of $C$
  have the requested property.\\
\end{remark}

\begin{example}
  The following figures show two examples of smooth tropical
  varieties:
  \begin{center}
    \includegraphics[scale=0.75]{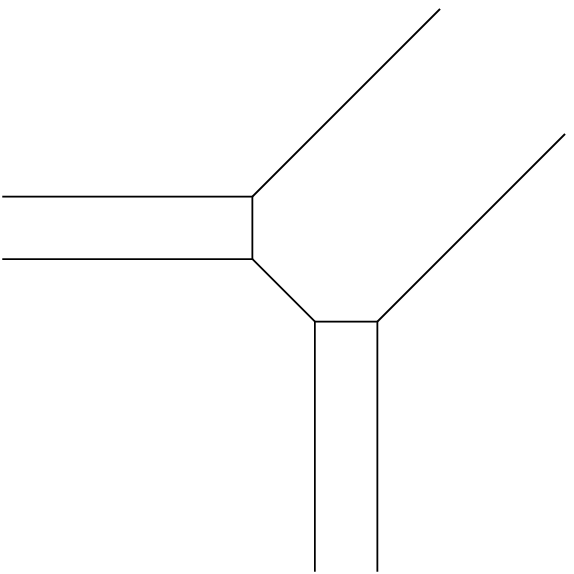}
    \qquad \qquad
    \includegraphics[scale=0.28]{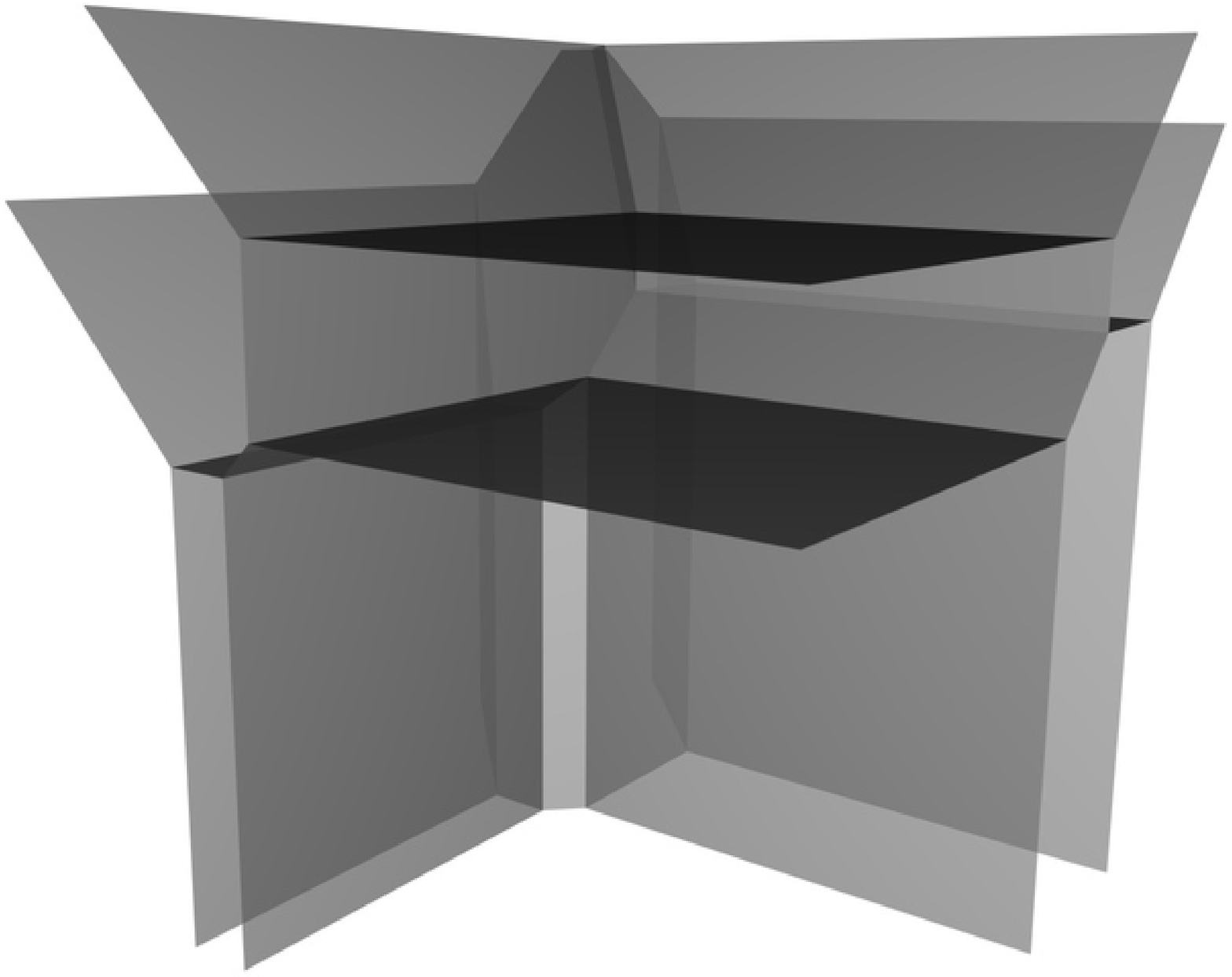}
  \end{center}
\end{example}

\begin{definition}
  Let $C$ be an abstract tropical cycle, $D$ a subcycle of $C$ with representative $X$
  and ${U \subseteq |C|}$ an open subset. We denote by $X \cap U$ the \emph{open tropical polyhedral
  complex} $$X \cap U := \left( \{ \sigma \cap U| \sigma \in X\}, |X| \cap U
  \right)$$ and by $\left[ X \cap U \right]$ its equivalence class modulo
  refinements. As this class only depends on the class of $X$ we
  can define $D \cap U :=\left[ X \cap U \right]$.
\end{definition}

\begin{remark}
  If we are given an open covering $\{U_1,\ldots,U_r\}$ of $C$ and
  open tropical polyhedral complexes $D_1 \cap U_1, \ldots, D_r \cap U_r$
  such that $D_i \cap U_i \cap U_j=D_j \cap U_i \cap U_j$ we can
  glue $D_1 \cap U_1, \ldots, D_r \cap U_r$ to obtain a cycle $D
  \in Z_{*}(C)$.
\end{remark}

\begin{definition}[Intersection products] \label{def-intersectionsmoothvariety}
  Let $C$ be a smooth tropical variety and let $(((X,|X|),\omega_X),\{ \Phi_\sigma\})$
  be a representative of $C$ like in definition \ref{def-smoothvarieties}.
  Moreover, let $D,E$ be subcycles of $C$. We construct
  local intersection products as follows: For every $\sigma \in X$
  we can regard $(D \cap S_\sigma)$ and $(E \cap S_\sigma)$ as open tropical cycles in
  $\widetilde{F_\sigma}$ via the map $\Phi_\sigma$. Let $\widetilde{D \cap S_\sigma}$ and $\widetilde{E \cap
  S_\sigma}$ be any tropical cycles in $\widetilde{F_\sigma}$ restricting to $D \cap S_\sigma$ and $E \cap S_\sigma$. As we have an
  intersection product on $\widetilde{F_\sigma}$ by remark \ref{remark-independence}
  we can define the intersection $$(D \cdot_\sigma E) \cap S_\sigma := \left( (\widetilde{D \cap S_\sigma}) \cdot (\widetilde{E \cap
  S_\sigma}) \right) \cap S_\sigma.$$
  Note that $(D \cdot_\sigma E) \cap S_\sigma$ does not depend on
  the choice of the cycles $\widetilde{D \cap S_\sigma}$ and $\widetilde{E \cap
  S_\sigma}$. Since $\{ S_\sigma | \sigma~ \in~ X \}$ is an open covering of
  $|C|$ and the local intersection products $(D \cdot_\sigma E) \cap
  S_\sigma$, $\sigma \in X$ are compatible by the following lemma we can
  glue them to obtain a global intersection cycle $D \cdot E \in Z_{*}(C)$.
\end{definition}

\begin{lemma} \label{lemma-gluingwelldefined}
  For the local intersection products in definition \ref{def-intersectionsmoothvariety} holds:
  $$(D \cdot_\sigma E) \cap S_\sigma \cap S_{\sigma'} = (D \cdot_{\sigma'} E) \cap S_\sigma \cap
  S_{\sigma'}.$$
\end{lemma}
\begin{proof}
  By definition we have an integer linear map $$|\widetilde{F_1}| \supseteq \Phi_\sigma(S_\sigma \cap
  S_{\sigma'}) \stackrel{f}{\longrightarrow} \Phi_{\sigma'}(S_\sigma \cap
  S_{\sigma'}) \subseteq |\widetilde{F_2}|$$
  with integer linear inverse $f^{-1}$,
  where $\widetilde{F_1}, \widetilde{F_2}$ are the tropical fans
  generated by \linebreak ${\Phi_\sigma(S_\sigma \cap S_{\sigma'})}$ and $\Phi_{\sigma'}(S_\sigma \cap
  S_{\sigma'})$, respectively. Let $C_1, C_2$ be subcycles of
  $\widetilde{F_1}$. We have to show that
  $$ C_1 \cdot C_2 = (f^{-1})_{*}(f_{*}(C_1) \cdot f_{*}(C_2)).$$
  If $\pi$ is the respective projection on the first factor we obtain by proposition
  \ref{propo-independence} and remark \ref{remark-independence} the equation
  $$\begin{array}{rcl}
    (f^{-1})_{*}(f_{*}(C_1) \cdot f_{*}(C_2))
    &=& (f^{-1})_{*} \left( \pi_{*} \left( \triangle_{\widetilde{F_2}} \cdot \left( f_{*}(C_1) \times f_{*}(C_2)\right) \right) \right)\\
    &=& \pi_{*} \left( (f^{-1} \times f^{-1})_{*} \left( \triangle_{\widetilde{F_2}} \cdot (f_{*}(C_1) \times f_{*}(C_2)) \right) \right)\\
    &=& \pi_{*} \left( (f^{-1} \times f^{-1})_{*} \left( (f \times f)_{*}(\triangle_{\widetilde{F_1}}) \cdot (f \times f)_{*}(C_1 \times C_2) \right) \right)\\
    &=& \pi_{*} \left( \triangle_{\widetilde{F_1}} \cdot C_1 \times C_2 \right)\\
    &=& C_1 \cdot C_2.
  \end{array}$$
\end{proof}

\begin{remark}
  Lemma \ref{lemma-gluingwelldefined} also implies that
  further refinements of the representative \linebreak $(((X,|X|),\omega_X),\{ \Phi_\sigma\})$
  of $C$ do not change the result $D \cdot E$. Hence the
  intersection product is well-defined.
\end{remark}

Our last step consists in proving basic properties of our
intersection product:

\begin{theorem}
  Let $C$ be a smooth tropical variety, let $D,D' \in Z_l(C)$, $E \in Z_{l'}(C)$ and
  $F \in Z_{l''}(C)$ be subcycles and let $\varphi \in \Div(C)$ be a Cartier
  divisor on $C$. Then the following equations hold in $Z_{*}(C)$:
  \begin{enumerate}
    \item $C \cdot D = D$,
    \item $D \cdot E = E \cdot D$,
    \item $(D+D') \cdot E = D \cdot E+ D' \cdot E,$
    \item $(D \cdot E) \cdot F = D \cdot (E \cdot F),$
    \item $\varphi \cdot (D \cdot E) = (\varphi \cdot D) \cdot E.$
  \end{enumerate}
  If moreover $D= (\sum_{i=1}^r \varphi_{i,1} \cdots \varphi_{i,l}) \cdot
  C$ for some Cartier divisors $\varphi_{i,j} \in \Div(C)$ then $$D
  \cdot E = \sum_{i=1}^r \varphi_{i,1} \cdots \varphi_{i,l} \cdot
  E$$ holds.
\end{theorem}
\begin{proof}
  The statements follow immediately from the definition of the intersection product
  and the corresponding statements in section \ref{sec-tropicallinearspaces}.
\end{proof}

\section{Pull-backs of cycles on smooth varieties} \label{sec-pullback}

We will now use the intersection product defined in section
\ref{sec-smoothvarieties} to introduce pull-backs of tropical
cycles along morphisms between smooth tropical varieties.

\begin{definition}[Pull-back] \label{def-smoothpullback}
  Let $X$ and $Y$ be smooth tropical varieties of dimension $m$ and $n$,
  respectively, and let $f:X \rightarrow Y$ be a morphism of tropical cycles.
  Moreover, let $\pi: X \times Y \rightarrow X$ be the projection onto the first factor
  and let $\gamma_f: X \rightarrow X \times Y$ be the morphism given by $x \mapsto (x,f(x))$.
  We denote by $\Gamma_f := (\gamma_f)_{*} X $ the graph of $f$.
  For a cycle $C \in Z_{n-k}(Y)$ we define its \emph{pull-back} $f^{*}C \in Z_{m-k}(X)$
  to be $$f^{*}C := \pi_{*} \left( \Gamma_f \cdot (X \times C) \right).$$
\end{definition}

The easiest non-trivial, but nevertheless important example of a
pull-back is the following:

\begin{example} \label{ex-pullbackwithprojection}
  Let $C$ and $D$ be smooth tropical cycles and let $p: C \times D \rightarrow D$
  be the projection on the second factor. We want to calculate the
  pull-back $p^{*}E$ for a cycle $E \in Z_k(D)$: The map $\gamma_p$
  from definition \ref{def-smoothpullback} is then just given by
  $\gamma_p: C \times D \rightarrow C \times D \times D: (x,y)
  \mapsto (x,y,y)$ and the map $\pi: C \times D \times D \rightarrow C \times D$
  is the projection to the first two factors. Hence we can conclude that
  $\Gamma_p=C \times \triangle_D$. Moreover, let
  $\pi^1: C \times D \times D \rightarrow C$ be the projection to the first and
  $\pi^2: C \times D \times D \rightarrow D$ be the projection to the
  second factor. We obtain by definition \ref{def-smoothpullback}:
  $$\begin{array}{rcl}
    p^{*}E &=& \pi_{*}(\Gamma_p \cdot (C \times D \times E))\\
    &=& \pi_{*}((C \times \triangle_D) \cdot (C \times D \times E))\\
    &=& \pi^1_{*}(C \cdot C) \times \pi^2_{*}(\triangle_D \cdot (D \times E))\\
    &=& C \times E.
  \end{array}$$
\end{example}

The pull-back has the following basic properties:

\begin{theorem} \label{thm-propertiesofpullback}
  Let $X, Y$ and $Z$ be smooth tropical varieties and let $f:X \rightarrow Y$ and $g:Y \rightarrow Z$ be morphisms of tropical cycles.
  Moreover, let $C, C' \in Z_{*}(Y)$ and $D \in Z_{*}(X)$ be
  subcycles. Then the following holds:
  \begin{enumerate}
    \item $f^{*}Y=X$,
    \item $\id_Y^{*}C=C$,
    \item if $C= \varphi_1 \cdots \varphi_r \cdot Y$ then $f^{*}C=f^{*} \varphi_1 \cdots f^{*} \varphi_r \cdot X$,
    \item $C \cdot f_{*}D = f_{*} (f^{*}C \cdot D)$,
    \item $(g \circ f)^{*} C= f^{*} g^{*} C$,
    \item $f^{*}(C \cdot C') = f^{*} C \cdot f^{*} C'$.
  \end{enumerate}
\end{theorem}
\begin{proof}
  Throughout the proof, let $\pi^X,\pi_X, \pi^{1}, \pi_{1}, \pi^Y,\pi_Y, \pi^{2}, \pi_{2}, \pi^{X,Y},\pi_{X,Y}, \pi^{1,2}, \pi_{1,2}$ and so forth be the projections to the
  respective factors.\\
  (a) and (b): By definition of the pull-back follows $$f^{*}Y = \pi^X_{*}(\Gamma_f \cdot (X \times Y))
  = \pi^X_{*}(\Gamma_f) = X$$ and $$\id_Y^{*}C = \pi^1_{*}(\Gamma_{\id_Y} \cdot (Y \times C))
  = \pi^1_{*}(\triangle_Y \cdot (Y \times C)) = Y \cdot C =C.$$
  (c): We have
  $$\begin{array}{rcl}
    f^{*}C &=& \pi_{*}^X \left( \Gamma_f \cdot (X \times (\varphi_1 \cdots \varphi_r \cdot Y)) \right)\\
    &=& \pi_{*}^X \left( \pi_2^{*}\varphi_1 \cdots \pi_2^{*}\varphi_r \cdot \Gamma_f \cdot (X \times Y) \right)\\
    &=& \pi_{*}^X \left( \pi_2^{*}\varphi_1 \cdots \pi_2^{*}\varphi_r \cdot \Gamma_f \right).
  \end{array}$$
  By definition of the intersection product (see \cite[definitions 3.4 and 6.5]{AR07})
  this last line is equal to $$f^{*} \varphi_1 \cdots f^{*} \varphi_r \cdot X.$$
  (d): Let $\pi_{X}: X \times Y \rightarrow X$ be the projection on
  $X$. By example \ref{ex-pullbackwithprojection} we know that $\pi_X^{*}D=D \times Y$. As the diagonal $\triangle_X$ can locally be expressed
  by Cartier divisors we can apply \cite[proposition 7.7]{AR07} and statement (c) locally to deduce
  that for all subcycles $E$ of $X \times Y$ holds
  $$\begin{array}{rcl}
    D \cdot \pi_{*}^X E &=& \pi_{*}^1 ( \triangle_X \cdot (D \times \pi_{*}^X E))\\
    &=& \pi_{*}^1 ( \triangle_X \cdot (\id \times \pi^X)_{*}(D \times E))\\
    &=& \pi_{*}^1 ((\id \times \pi^X)_{*}( (\id \times \pi^X)^{*}\triangle_X \cdot (D \times E)))\\
    &=& \pi_{*}^1 ((\id \times \pi^X)_{*}( (\triangle_X \times Y) \cdot (D \times E)))\\
    &=& \pi_{*}^1 (\pi_{*}^{1,2} (  (\triangle_X \times Y) \cdot (D \times E)))\\
    &=& \pi_{*}^1 (\pi_{*}^{1,2} (  \triangle_{X \times Y} \cdot (D \times Y \times E)))\\
    &=& \pi_{*}^1 ((D \times Y) \cdot E)\\
    &=& \pi_{*}^X ( \pi^{*}_X D \cdot E).
  \end{array}$$
  This implies that
  $$\begin{array}{rcl}
    f^{*}C \cdot D &=& D \cdot \pi_{*}^X(\Gamma_f \cdot (X \times C))\\
    &=& \pi_{*}^X( \pi^{*}_X D \cdot \Gamma_f \cdot (X \times C))\\
    &=& \pi_{*}^X( (D \times Y) \cdot \Gamma_f \cdot (X \times C))\\
    &=& \pi_{*}^X( \Gamma_f \cdot (D \times C)).
  \end{array}$$
  Moreover, it is easy to check that $(f \times \id)^{*} \triangle_Y =
  \Gamma_f$. As above we can conclude that
  $$\begin{array}{rcl}
    C \cdot f_{*} D &=& \pi_{*}^1 ( \triangle_Y \cdot (C \times f_{*} D))\\
    &=& \pi_{*}^1 ( (\id \times f)_{*}( (\id \times f)^{*}\triangle_Y \cdot (C \times D) ))\\
    &=& f_{*} ( \pi_{*}^X ( (\id \times f)^{*}\triangle_Y \cdot (C \times D) ))\\
    &=& f_{*} ( \pi_{*}^X ( (f \times \id)^{*}\triangle_Y \cdot (D \times C) ))\\
    &=& f_{*} ( \pi_{*}^X ( \Gamma_f \cdot (D \times C) ))\\
    &=& f_{*} ( f^{*}C \cdot D ).
  \end{array}$$
  (e): Let $\Phi: X \rightarrow X \times Y \times Z$ be given by
  $x \mapsto (x,f(x),g(f(x)))$. An easy calculation shows that
  $(\Gamma_f \times Z) \cdot (X \times \Gamma_g) = \Phi_{*}X$. Hence we can conclude by statement (d) that
  $$\begin{array}{rcl}
    f^{*} g^{*} C &=& \pi_{*}^X \left( \Gamma_f \cdot \left(X \times \pi_{*}^Y \left( \Gamma_g \cdot (Y \times C) \right) \right) \right)\\
    &=& \pi_{*}^X \left( \pi_{*}^{X,Y} \left( (\Gamma_f \times Z) \cdot (X \times \Gamma_g) \cdot (X \times Y \times C) \right) \right)\\
    &=& \pi_{*}^X \left( (\Gamma_f \times Z) \cdot (X \times \Gamma_g) \cdot (X \times Y \times C) \right)\\
    &=& \pi_{*}^X \left( \Phi_{*}X \cdot (X \times Y \times C) \right)\\
    &=& \pi_{*}^X \left( \Gamma_{g \circ f} \cdot (X \times C) \right)\\
    &=& (g \circ f)^{*} C.
  \end{array}$$
  (f): Let $\Phi:X \rightarrow X \times Y \times Y$ be given by $x
  \mapsto (x,f(x),f(x))$ and let $\pi^{1,2}, \pi^{1,3}: X \times Y \times Y \rightarrow X \times Y$ be the
  projections to the respective factors. An easy calculation shows that $$(\Gamma_f
  \times Y) \cdot (X \times \Gamma_{\id_Y}) = \Phi_{*} X = \pi_{1,2}^{*} \Gamma_f \cdot \pi_{1,3}^{*}
  \Gamma_f.$$ Hence we can deduce that
  $$\begin{array}{rcl}
    f^{*}(C \cdot C') &=& \pi_{*}^X \left( \Gamma_f \cdot (X \times (C \cdot C') ) \right)\\
    &=& \pi_{*}^X \left( \Gamma_f \cdot (X \times \pi_{*}^1 (\Gamma_{\id_Y} \cdot C \times C') ) \right)\\
    &=& \pi_{*}^X \left( \Gamma_f \cdot \pi_{*}^{1,2} ( (X \times \Gamma_{\id_Y}) \cdot (X \times C \times C') ) \right)\\
    &=& \pi_{*}^X \left( \pi_{*}^{1,2} ( (\Gamma_f \times Y) \cdot (X \times \Gamma_{\id_Y}) \cdot (X \times C \times C') ) \right)\\
    &=& \pi_{*}^X \left( \pi_{*}^{1,3} ( (\Gamma_f \times Y) \cdot (X \times \Gamma_{\id_Y}) \cdot (X \times C \times C') ) \right)\\
    &=& \pi_{*}^X \left( \pi_{*}^{1,3} (  \pi_{1,2}^{*} \Gamma_f \cdot \pi_{1,3}^{*} \Gamma_f \cdot (X \times C \times C') ) \right)\\
    &=& \pi_{*}^X \left( \Gamma_f \cdot \pi_{*}^{1,3} ( (\Gamma_f \times Y) \cdot (X \times C \times C') ) \right)\\
    &=& \pi_{*}^X \left( \Gamma_f \cdot (\pi_{*}^{X} ( \Gamma_f \cdot (X \times C)) \times C') \right)\\
    &=& \pi_{*}^X \left( \Gamma_f \cdot (f^{*}C \times C') \right)\\
    &=& f^{*}C \cdot f^{*} C'.
  \end{array}$$
\end{proof}

We finish the section with another important example:

\begin{example}
  Let $D$ be a smooth tropical variety and let $C \in Z_k(D)$ be a smooth tropical subvariety.
  Moreover, let $\iota: C \rightarrow D$ be the inclusion map. We want to calculate
  the pull-back $\iota^{*}E$ for a cycle $E \in Z_l(D)$: Let
  $\pi^C:C \times D \rightarrow C$ and $\pi^{D}:C \times D \rightarrow D$ be
  the projections to the first and second factor and let
  $\gamma_\iota: C \rightarrow C \times D$ be given by $x \mapsto (x,x)$.
  Hence we can deduce that $\Gamma_\iota=(\gamma_\iota)_{*}C=\triangle_C$
  and by example \ref{ex-pullbackwithprojection} that
  $(\pi^{D})^{*}E=C \times E$.
  Thus we can conclude by theorem \ref{thm-propertiesofpullback} (d):
  $$\begin{array}{rcl}
    \iota^{*}E &=& \pi^C_{*}(\Gamma_\iota \cdot (C \times E))\\
    &=& \pi^C_{*}(\triangle_C \cdot (C \times E))\\
    &=& \pi^{D}_{*}(\triangle_C \cdot (C \times E))\\
    &=& \pi^{D}_{*}(\triangle_C \cdot (\pi^{D})^{*}E)\\
    &=& \pi^{D}_{*}(\triangle_C) \cdot E\\
    &=& C \cdot E,
  \end{array}$$
  where $C \cdot E$ is the intersection product on $D$.
\end{example}

I would like to thank my advisor Andreas Gathmann for numerous
helpful discussions.


\begin {thebibliography}{MMMM}

\bibitem [AR07]{AR07}
  \arxivjournal{Lars Allermann, Johannes Rau}
               {First steps in tropical intersection theory}
               {Mathematische Zeitschrift (to appear)}
               {0709.3705}

\bibitem [AR08]{AR08}
  \arxiv{Lars Allermann, Johannes Rau}
        {Tropical rational equivalence}
        {0811.2860}

\bibitem [GKM07]{GKM07}
  \arxivjournal{Andreas Gathmann, Michael Kerber, Hannah Markwig}
               {Tropical fans and the moduli spaces of tropical curves}
               {Compositio Mathematica, Volume 145, 2009, no. 1, pp. 173--195(23)}
               {0708.2268}

\bibitem [M06]{M06}
  \arxivjournal{Grigory Mikhalkin}
               {Tropical Geometry and its applications}
               {Proceedings of the ICM, Madrid, Spain, 2006, pp. 827--852(26)}
               {math.AG/0601041}

\bibitem [M07]{M07}
  \arxiv{Grigory Mikhalkin}
        {Introduction to Tropical Geometry, notes from the IMPA lectures, summer 2007}
        {0709.1049}

\end {thebibliography}
\end {document}